\documentclass{article}
\usepackage{amssymb}

\newtheorem{thm}{Theorem}

\newtheorem{lem}[thm]{Lemma}

\newenvironment{proof}{
{\vspace{.3cm}\noindent\sc{Proof:}}\
}{\par}

\renewcommand{\(}{\left(}
\renewcommand{\)}{\right)}

\newcommand{\got}{\mathfrak}

\newcommand{\X}{X_1,\ldots ,X_n}

\newcommand{\Xm}{X_1,\ldots ,X_m}

\renewcommand{\v}{\wh{v}}

\newcommand{\wh}{\widehat}

\newcommand{\al}{\alpha}

\newcommand{\R}{R_{\v}}
\renewcommand{\r}{R_{v}}
\newcommand{\M}{{\got m}_{\v}}
\newcommand{\m}{{\got m}_v}

\newcommand{\Kn}{{\widehat K_n}}
\newcommand{\D}{\Delta_{\v}}
\renewcommand{\d}{\Delta_v}

\newcommand{\lcor}{{\rm [\kern - 1.8pt [}}
\newcommand{\rcor}{{\rm ]\kern - 1.8pt ]}}
\renewcommand{\[}{\left[}
\renewcommand{\]}{\right]}
\newcommand{\llcor}{\[\kern -2.5pt\[}
\newcommand{\rrcor}{\]\kern -2.5pt\]}

\newcommand{\fb}{\overline{f}}

\newcommand{\ZZ}{{\mathbb{Z}}}

\def\boxit [#1]#2{\vbox{\hrule\hbox{\vrule
     \vbox spread #1{\vss\hbox spread #1{\hss #2\hss}\vss}%
      \vrule}\hrule}}

\def\qed{\boxit[5pt]{}} 

\begin{document}

\vspace{2.5cm}

\title{ON THE DIMENSION OF DISCRETE VALUATIONS OF $k((\X ))$}

\author{ Miguel \'{A}ngel Olalla Acosta 
\thanks{Partially supported by Junta de Andaluc{\'{\i}}a, Ayuda a
grupos FQM 218.} 
}





\date{Departamento de \'{A}lgebra. Universidad de Sevilla, Espa\~{n}a.}



\maketitle

\vspace{1.5cm}


\begin{abstract}
Let $v$ be a rank-one discrete valuation of the field $k((\X ))$. We
know, after \cite{Bri2}, that if $n=2$ then  the dimension
of $v$ is 1
and if $v$ is the usual order function over $k((\X ))$ its dimension
is $n-1$. In this paper we prove that, in the general case, the
dimension of a rank-one discrete valuation can be any number between 1
and $n-1$.
\end{abstract}

\thispagestyle{empty}

\section{TERMINOLOGY AND PRELIMINARIES}

Let $k$ be an algebraically closed field of characteristic 0,
$R_n=k\lcor\X\rcor$, $M_n=(\X )$ the maximal ideal and
$K_n=k ((\X ))$ the quotient field. Let $v$ be a rank-one discrete
valuation of $K_n\vert k$, $\r$ its valuation ring, $\m$ its
maximal ideal and $\d$ its residual field of $v$. The center of
$v$ in $R_n$ is $\m\cap R_n$. Throughout this paper ``discrete
valuation of $K_n\vert k$" will mean ``rank-one discrete valuation of
$K_n\vert k$ whose center in $R_n$ be the maximal ideal $M_n$".
The dimension of $v$ is the transcendence degree of $\d$ over $k$.
We shall suppose that the group of $v$ is $\ZZ$.

Let $\Kn$ be the completion of $K_n$ with respect to $v$ (see \cite{Ser1}), $\v$ the
extension of $v$ to $\Kn$, $\R$, $\M$ and $\D$ the ring, maximal
ideal and the residual field of $\v$, respectively. 

\section{THE DIMENSION OF $v$}

Fix a number $m$ between 1 and $n-1$. We are featuring a constructive
method in order to obtain examples of valuations with dimension $m$.

Let us consider the following homomorphism:
$$\begin{array}{crcl}
  \varphi : & k\lcor \X\rcor & \longrightarrow & \overline{k
                                                 (u)}\lcor t\rcor\\
      &                 X_1 & \longmapsto & t \\
      &                 X_2 & \longmapsto & ut \\
      &                 X_i & \longmapsto &\sum_{j\geq 1}u^{1/p_i^j}t^j
  \end{array}
$$
where $\overline{k (u)}$ stands for the algebraic closure of $k(u)$ and
$2<p_3<\ldots <p_n$ are prime numbers.

\begin{lem}
The homomorphism $\varphi$ is one to one.
\end{lem}

\begin{proof}
Let us take the fields $K_2=k(u)$ and $$K_i=k
(u,\{u^{1/p_3^j},j\geq 1\} ,\ldots ,\{u^{1/p_i^j},j\geq 1\})$$ for
all $i\geq 3$.

Let us suppose that $\varphi$ is not one to one, then $\ker
(\varphi )\ne\{ 0\}$. So let $f$ be a non-zero element of $M=(\X )$ such
that $f\in\ker (\varphi )$. Let $m$ be the higher index such that $f\in
k\lcor\Xm\rcor$.

\vspace{0.3cm}

\noindent
If $\underline{m=1\mbox{ {or} }2}$, trivially we have a contradiction.

\vspace{0.3cm}

\noindent
If $\underline{m=3}$, let us take
$$\fb =f(\varphi (X_1),\varphi (X_2),X_3)\in K_2\lcor t, X_3\rcor ,$$
and consider the homomorphism
$$\begin{array}{crcl}
  \psi : & K_2\lcor t,X_3\rcor & \longrightarrow & \overline{k
                                                 (u)}\lcor t\rcor\\
      &                 t & \longmapsto & t \\
      &                 X_3 & \longmapsto &\sum_{j\geq 1}u^{1/p_3^j}t^j.
  \end{array}
$$ We know that $\fb\in\ker (\psi )$ and this kernel is a prime ideal
because $\psi$ is an homomorphism between integral domains. We
can write $\fb =t^rg$, with $r\geq 0$ and $t$ doesn't divide to
$g$. This forces $g$ to have some non-trivial terms in $X_3$. Let
$s>0$ be the minimum such that $\al X_3^s$ is one of these terms. By the
Weierstrass preparation theorem we have $g=Ug'$, where $U(t,X_3)$
is a unit and $$g'=X_3^s+a_1(t)X_3^{s-1}+\ldots +a_s(t).$$ Since
$U$ is a unit, $g'\in\ker (\psi )$ and $$\psi
(g')=g'\left(t,\sum_{j\geq 1}u^{1/p_3^j}t^j\right)=0.$$ 
This leads to a contradiction because the roots of $g'$ are in $K_2\lcor
t^{1/q}\rcor$, with $q\in \ZZ$, by the Puiseux theorem.

\vspace{0.3cm}

\noindent
If \underline{$m>3$} let us take
$$\overline{f}=f(\varphi (X_1),\ldots ,\varphi (X_{m-1}),X_m)\in
K_{m-1}\lcor t,X_m\rcor$$
and consider the homomorphism
$$\begin{array}{crcl}
  \psi : & K_{m-1}\lcor t,X_m\rcor & \longrightarrow & \overline{k
                                                 (u)}\lcor t\rcor\\
      &                 t & \longmapsto & t \\
      &                 X_m & \longmapsto &\sum_{j\geq 1}u^{1/p_m^j}t^j.
  \end{array}
$$ 
As in the previous case we can write $\fb =t^rh$,
where $h\in\ker (\psi )$. So we have $h=Uh'$, where $U(t,X_m)$ is
a unit and $$h'=X_m^r+b_1(t)X_m^{r-1}+\ldots +b_r(t)\in \ker
(\psi ),$$ so $$h'\left(t,\sum_{j\geq 1}u^{1/p_m^j}t^j\right)=0.$$
But this is again a contradiction by the Puiseux theorem: since $\ker
(\psi )$ is a prime ideal, we can suppose that $h'$ is an
irreducible element of the ring $K_{m-1}\lcor t\rcor [X_m]$. In
this situation the Puiseux theorem says that to obtain the coefficients of
a root of $h'=0$, like a Puiseux series in $t$ with coefficients
in $k (u)$, we have to resolve {\em a finite number} of algebraic
equations of degree greater than 1 in $K_{m-1}$. Inside $K_{m-1}$
we can not obtain $u^{1/p_m}$ and, with a finite number of
algebraic equations, we can obtain {\em a finite number} of powers
of $u^{1/p_m^j}$ but not all. So this proves the lemma.\hfill\qed
\end{proof}

We shall extend to the quotient fields this injective homomorphism for
giving an example of a rank-one discrete valuation of $k((\X
))$ of dimension 1.

\begin{lem}
There exists a rank-one discrete valuation of $k ((\X ))$ of dimension 1.
\end{lem}

\begin{proof}
We know that the homomorphism $$\varphi :k \lcor\X\rcor\to\overline{k
(u)}\lcor t\rcor$$ previously defined is one to one. So we can take the
valuation $v=\nu\circ\varphi$, where $\nu$ is the usual order
function over $\overline{K(u)}((t))$ in $t$ and $\varphi$ is the
natural extension to the quotient fields. Let $\al$ be the
residue $X_2/X_1+\m\in\d$. Hence, to obtain the lemma we have to prove
that $\al\notin k$ and $\d$ is an algebraic extension of $k(\al
)$.

Let us suppose that $\al\in k$. Then there must exist $a\in k$ such that
$X_2/X_1+\m = a+\m$, so $$\frac{X_2-aX_1}{X_1}\in{\got m}_v.$$
This means that $v(X_2-aX_1)>1$. On the other side we have
$$\varphi (X_2-aX_1)=(u-a)t,$$ so $v(X_2-aX_1)=1$ and we have a
contradiction. Hence $\al\notin k$.

Let us prove that $\d$ is an algebraic extension of $k(\al )$. We
can consider each element of $k\lcor\X\rcor$ like a sum of forms
with respect to the usual degree. If $f_r$ is a form of degree $r$, 
then $\varphi (f_r)=t^rP$, with $P$ a polynomial in $u$ and a finite
number of elements $u^{1/p_i}$.

Let us take $f,\ g\in k\lcor\X\rcor$ such that $g\ne 0$ and
$v(f/g)=0$. Then $\varphi (f/g)=h_0+th_1$, where $h_0$ is a
rational fraction in $u$ and a finite number of elements
$u^{1/p_i}$. So $h_0$ is algebraic over $k(u)$. Let us consider
$$P(u,Z)=c_0(u)Z^m+c_1(u)Z^{m-1}+\ldots +c_{m-1}(u)Z+c_m(u)\in k[u][Z]$$ 
a polynomial satisfied by $h_0$, where $c_i(u)\in k[u]$ for all $i$ and
$c_0\ne 0$. Let $\beta$ be the element 
$$ \beta
=P\left(\frac{X_2}{X_1},\frac{f}{g}\right)
         =c_0\left(\frac{X_2}{X_1}\right)\left(\frac{f}{g}\right)^m+
          \ldots +c_m\left(\frac{X_2}{X_1}\right) .
$$
Then we have
$$ \varphi (\beta )=c_0(u)(h_0+th_1)^m+\ldots +c_m(u), $$
so $v(\beta )=\nu\circ\varphi (\beta) >0$ and $\beta\in\m$. 
Subsequently,
$$ 0+{\got m}_v=\beta+{\got m}_v=P\left(\alpha ,\frac{f}{g}+{\got
                                 m}_v\right) .$$
This proves that $f/g+\m$ is an algebraic element over $k(\al )$ and,
a fortiori, the lemma.\hfill\qed
\end{proof}

\begin{lem}
The dimension of a rank-one discrete valuation of $k((\X
))$ is between 1 and $n-1$.
\end{lem}

\begin{proof}
We know, after \cite{Bri2}, that the dimension of a
rank-one discrete
valuation of $k((\X ))$ is minor or equal than $n-1$. So we have
to prove that there exists a transcendental residue in $\d$.

Let us suppose that $v(X_i)=n_i$ for all $i=1,\ldots ,n$. Then the
value of $X_2^{n_1}/X_1^{n_2}$ is zero, so $0\ne
(X_2^{n_1}/X_1^{n_2})+\m\in\d$. If this residue lies in $k$ then
there exists $a_{21}\in k$ such that $$\frac{X_2^{n_1}}{X_1^{n_2}}+\m =
a_{21}+\m .$$ 
This implies
$$\frac{X_2^{n_1}}{X_1^{n_2}}-a_{21}=\frac{X_2^{n_1}-a_{21}X_1^{n_2}}{X_1^{n_2}}
\in\m ,$$ and then
$$v\(\frac{X_2^{n_1}-a_{21}X_1^{n_2}}{X_1^{n_2}}\) >0.$$ So we
have $v(X_2^{n_1}-a_{21}X_1^{n_2})=m_1>n_1n_2$. Then
$$v\(\frac{(X_2^{n_1}-a_{21}X_1^{n_2})^{n_1}}{X_1^{m_1}}\) =0.$$
If the residue of this element lies too in $k$, then there must exist
$a_{22}\in k$ such that
$v((X_2^{n_1}-a_{21}X_1^{n_2})^{n_1}-a_{22}X_1^{m_1})=m_2>n_1m_1$.
We can repeat this operation.

The previous procedure is finite: if it didn't stop we would construct
the power series 
$$X_2^{n_1}-\sum_{i=1}^{\infty}b_{2i}X_1^{i}$$ 
such that the sequence of partial sums has increasing values. Since
$\wh{K_n}$ is a complete field, then this series amounts to zero in
contradiction with $X_1$ and $X_2$ being formally independent. So the
procedure must stop and there exists a transcendental element over $k$ in $\d$.
\hfill\qed
\end{proof}

\begin{thm}
Let $m$ be a fixed number between 1 and $n-1$, then there exists a
rank-one discrete valuation of $k((\X ))$ of dimension $m$.
\end{thm}

\begin{proof}
Let us consider the one to one (the proof of injectivity
parallels that of lemma 1) homomorphism
$$\begin{array}{crcl}
\varphi & k\lcor\X\rcor & \longrightarrow
    & \overline{k (u)}\lcor t_1,\ldots ,t_m\rcor \\
 & X_1 & \longmapsto & t_1 \\
 & X_2 & \longmapsto & ut_1 \\
 & X_i & \longmapsto &
   \left\{\begin{array}{l}
      t_{i}\mbox{ if } i\leq m+1 \\
      \sum_{j\geq 1} u^{1/p_i^j}t_1^j\mbox{ if } i> m+1.
         \end{array}\right.
\end{array}
$$ We can take the valuation $v:=\nu\circ\varphi$, with $\nu$ the
usual order function in $\overline{k(u)}\lcor t_1,\ldots
,t_m\rcor$ and $\varphi$ the natural extension to the quotient
fields. We know (lemma 2) that the residue $X_2/X_1+\m$ is
transcendental over $k$. Trivially the residue $X_i/X_1+\m$ for
all $i=3,\ldots ,m+1$ are transcendental over $k(X_2/X_1+\m ,\ldots
,X_{i-1}+\m )$ because $t_i$ are formally independent variables.
Any element $f/g+\m\in\d$ is algebraic over $k(X_2/X_1+\m ,\ldots
,X_{m+1}+\m )$ parallels that of lemma 2. So the dimension of $v$ is $m$.
\hfill\qed
\end{proof}



\providecommand{\bysame}{\leavevmode\hbox to3em{\hrulefill}\thinspace}

\end{document}